\documentclass[11pt]{article}
\usepackage{amssymb,amsfonts,amsmath,latexsym,epsf,tikz,url}

\newtheorem{theorem}{Theorem}[section]
\newtheorem{proposition}[theorem]{Proposition}

\newtheorem{conjecture}[theorem]{Conjecture}
\newtheorem{corollary}[theorem]{Corollary}
\newtheorem{lemma}[theorem]{Lemma}

\newtheorem{definition}[theorem]{Definition}

\newcommand{\proof}{\noindent{\bf Proof.\ }}
\newcommand{\qed}{\hfill $\square$\medskip}

\begin{document}

\title{Stabilizing on the distinguishing number of a graph}

\author{
Saeid Alikhani  $^{}$\footnote{Corresponding author}
\and
Samaneh Soltani
}

\date{\today}

\maketitle

\begin{center}
Department of Mathematics, Yazd University, 89195-741, Yazd, Iran\\
{\tt alikhani@yazd.ac.ir, s.soltani1979@gmail.com}
\end{center}


\begin{abstract}
The distinguishing number $D(G)$  of a graph $G$ is the least integer $d$
such that $G$ has a vertex labeling   with $d$ labels  that is preserved only by a trivial
automorphism. The distinguishing  stability, 
 of a graph $G$ is denoted by $st_D(G)$ and is the minimum number of vertices whose
removal changes the distinguishing  number.  We obtain  a general upper bound $st_D(G) \leqslant \vert V(G)\vert -D(G)+1$, and a  relationships between the distinguishing stabilities of graphs $G$ and $G-v$, i.e., $st_D(G)\leqslant st_D(G-v)+1$, where $v\in V(G)$. Also we study  the edge distinguishing stability number (distinguishing bondage number) of $G$. 
\end{abstract}


\medskip
\noindent{\bf AMS Subj.\ Class.}: 05C25 

\section{Introduction and definitions}

Let $G=(V,E)$ be a simple graph of order $n\geqslant 2$. We use the the following notations: The set of vertices adjacent in $G$ to a vertex of a vertex subset $W\subseteq  V$ is the open neighborhood $N_G(W )$ of $W$. The closed neighborhood  $G[W ]$ also includes all vertices of $W$ itself. In case of a singleton set $W =\{v\}$ we write $N_G(v)$ and $N_G[v]$ instead of $N_G(\{v\})$ and $N_G[\{v\}]$, respectively. ${\rm Aut}(G)$ denotes the automorphism group of $G$.  
A labeling of $G$, $\phi : V \rightarrow \{1, 2, \ldots , r\}$, is said to be $r$-distinguishing, 
if no non-trivial  automorphism of $G$ preserves all of the vertex labels.
The point of the labels on the vertices is to destroy the symmetries of the
graph, that is, to make the automorphism group of the labeled graph trivial.
Formally, $\phi$ is $r$-distinguishing if for every non-trivial $\sigma \in {\rm Aut}(G)$, there
exists $x$ in $V$ such that $\phi(x) \neq \phi(x\sigma)$. The distinguishing number of a graph $G$ is defined  by
\begin{equation*}
D(G) = {\rm min}\{r \vert ~ G ~\text{{\rm has a labeling that is $r$-distinguishing}}\}.
\end{equation*} 

This number has defined by Albertson and Collins \cite{Albert}. If a graph has no nontrivial automorphisms, its distinguishing number is one. In other words, $D(G) = 1$ for the asymmetric graphs. The other extreme, $D(G) = \vert V(G) \vert$, occurs if and only if $G = K_n$. The distinguishing number of some examples of graphs, $D(P_n) =2$ for every $n\geqslant 3$, and  $D(C_n) =3$ for $n =3,4,5$,  $D(C_n) = 2$ for $n \geqslant 6$. Also $D(K_{p,q})=p$, for $p > q$, and $D(K_{p,p})=p+1$, for $p\geq 4$. 
Authors in \cite{operation} have shown that removing a vertex of $G$ can decrease the distinguishing number by at most one but can increase by at most to double  of distinguishing number  of $G$.  
Also for each connected graph  $G$ and $e\in E(G)$, $\vert D(G-e)-D(G) \vert \leqslant 2$.

A domination-critical (domination-super critical, respectively) vertex in a graph
$G$ is a vertex whose removal decreases (increases, respectively) the domination
number. Bauer et al. \cite{1}  introduced the concept of domination stability in graphs.
The domination stability, or just $\gamma$-stability, of a graph $G$ is the minimum number
of vertices whose removal changes the domination number. Motivated by domination stability, we introduce the distinguishing stability of a graph. 

\begin{definition}
	Let $G$ be a graph of order $n\geqslant 2$. The stabilizing on the distinguishing number, or just distinguishing stability,  $st_D(G)$ of graph $G$ is the  minimum number of vertices whose
	removal changes the distinguishing  number.
\end{definition}
 
Also we introduce and study  the edge distinguishing stability number (distinguishing bondage number) of $G$ and compute edge distinguishing stability of some specific graphs. 

In the next section we compute the distinguishing stability of some specific graphs.   We obtain general bounds, and a relationships between the distinguishing stabilities   of $G$ and $G-v$, where $G-v$ denotes the graph obtained from $G$ by removal of a vertex  $v$ and all edges incident to $v$, in Section 3. Finally we consider and study the edge distinguishing stability number of graphs in Section 4. 

\section{Distinguishing stability of specific graphs} 

\begin{figure}[ht]
	\begin{center}
		\includegraphics[width=0.4\textwidth]{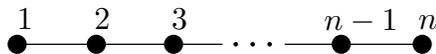}
		\caption{\label{fpath} {\small A path $P_n$ with vertex set $\{1,...,n\}$.}}
	\end{center} 
\end{figure}

In this section, first we compute the distinguishing stability of some specific graphs. We start with paths and cycles. A path is a connected
graph in which two vertices have degree one  and the remaining vertices have degree two. Let
$P_n$ be the path with $n$ vertices as shown in Figure \ref{fpath}.

\begin{proposition}\label{path}
For any $n\geqslant 6$, the distinguishing stability of $P_n$, is  $st_D(P_n)=2$, while $st_D(P_2)=1$, $st_D(P_3)=2$, $st_D(P_4)=3$, and $st_D(P_5)=1$.
\end{proposition}
\proof  It is clear that $st_D(P_n) > 1$ for $n\geqslant 6$. On the other hand by removing the third and sixth vertex of $P_n$, i.e., vertices labeled by numbers $3$ and $6$, the graph $P_n$ convert to disjoint union of the two paths $P_2$ and one path $P_{n-5}$, and so $D(P_n-\{3,6\})=3$. Since $D(P_n) =2$, so $st_D(P_n)=2$ for $n\geqslant 6$.\qed

The following proposition obtain immediately from Proposition \ref{path}.
\begin{proposition}
For any $n\geqslant 7$, the distinguishing stability of $C_n$ is, $st_D(C_n)=3$, while $st_D(C_3)= st_D(C_4)= st_D(C_5)= 1$, and $st_D(C_6)=2$.
\end{proposition}

With respect to the value of the distinguishing number of the complete graphs and complete bipartite graphs we can prove the following proposition:
\begin{proposition}
	\begin{enumerate}
		\item[(i)] For any $p\geqslant 1$, $st_D(K_p)=1$.

\item[(ii)]   For any $n\geqslant m$, $st_D(K_{n,m})=1$, except, $ st_D(K_{n,n+1})=2$.
\end{enumerate} 
\end{proposition}

\begin{figure}
	\begin{center}
		\hspace{.7cm}
		\includegraphics[width=0.5\textwidth]{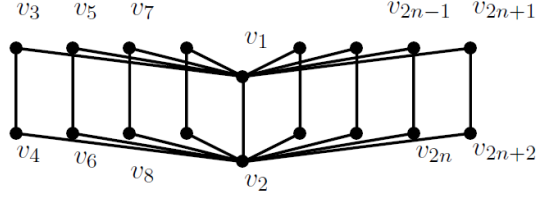}
				\caption{\label{book}  Book graph $B_n$.}
	\end{center}
\end{figure}

 The $n$-book graph $(n\geqslant 2)$ (Figure \ref{book}) is defined as the Cartesian product $K_{1,n}\square P_2$. We call every $C_4$ in the book graph $B_n$, a page of $B_n$. All pages in $B_n$ have a common side $v_1v_2$.   We shall compute the distinguishing stability number  of $B_n$. 
 The following theorem gives the distinguishing number of the book graph. 
 
\begin{theorem}{\rm \cite{Alikhani}}
	The distinguishing number of $B_n$ $(n\geq 2)$ is $D(B_n)= \lceil \sqrt{n}\rceil $.
\end{theorem}

\begin{proposition}
The distinguishing stability of the book graph is 
\begin{equation*}
 st_D(B_n)= \left\{\begin{array}{ll}
1& \text{{\rm if $n-1$  is square}},\\
2 & {\rm otherwise}.
\end{array}\right.
\end{equation*}
\end{proposition}
\proof  By removing the two central vertices of $B_n$, the book graph $B_n$ convert to disjoint union of $n$ paths $P_2$, denoted by $nP_2$. It can be computed that the distinguishing number of $nP_2$ is $D(nP_2)= \lceil \dfrac{1+\sqrt{8n+1}}{2}\rceil$ which  is different from the distinguishing number of $B_n$, $D(B_n)= \lceil \sqrt{n}\rceil $. So $st_D(B_n)\leqslant 2$. If the value of $n-1$ is square, then $D(B_n)= D(B_{n-1})+1$, and so by removing a noncentral vertex of $B_n$, say $v$, we have a book graph $B_{n-1}$ and the path $P_3$, that start point of $P_3$ is identified with one of the central point of $B_{n-1}$, and so $P_3$ is fixed by any automorphism of $B_n - v$. So $D(B_n - v)= D(B_{n-1})$, and hence $st_D(B_n)=1$.\qed

\begin{figure}[ht]
	\begin{center}
		\hspace{.7cm}
		\includegraphics[width=0.5\textwidth]{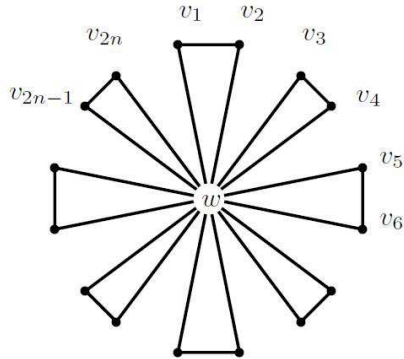}
				\caption{\label{friendsep}  Friendship graph $F_n$.}
	\end{center}
\end{figure}

The friendship graph $F_n$ $(n\geqslant 2)$ can be constructed by joining $n$ copies of the cycle graph $C_3$ with a common vertex. To compute the distinguishing stability of $F_n$, we need the following theorem.  
\begin{theorem}{\rm \cite{Alikhani}}\label{distfan}
The distinguishing number of the friendship graph $F_n$  $(n\geq 2)$ is  $$D(F_n)= \lceil \dfrac{1+\sqrt{8n+1}}{2}\rceil.$$  
\end{theorem}
\begin{proposition}\label{pro fan}
The distinguishing stability of the friendship graph is  $$st_D(F_n)= {\rm min}\{k: 8(n-k)+1~{\rm is~square}\}.$$
\end{proposition}
\proof  Set ${\rm min}\{k: 8(n-k)+1~{\rm is~square}\}=t$ and consider the graph $F_n$ as shown in Figure \ref{friendsep}. The graph  
$F_n-\{v_1,v_3,\ldots ,v_{2t-1}\}$ is two graphs $F_{n-t}$ and $K_{1,t}$ such that their central vertices are identified. Since $t$ is the minimum number which $ 8(n-t)+1$ is square, so $D(F_{n-t})=D(F_n)-1$. On the other hand, since $K_{1,t}$ and $F_{n-t}$ are two nonisomorphic graphs, so $D(F_n-v_1- v_3-\cdots - v_{2t-1})= {\rm max}\{D(F_{n-t}) , D(K_{1,t})\}$, hence $D(F_n-v_1- v_3-\cdots - v_{2t-1}) =D(F_n)-1$. Thus $st_D(F_n)\leqslant t$. Now we show that $st_D(F_n) >  t-1$. First note that if we remove the central vertex of $F_n$, say $w$, then the value of the distinguishing number of $F_n -w$ is equal $D(F_n)$. In general, if $v_{i_1},\ldots , v_{i_s}$ are $s$ noncentral vertices of $F_n$, then $D(F_n - v_{i_1}- \cdots - v_{i_s})= D(F_n - v_{i_1}- \cdots - v_{i_s}-w)$, because the central vertex $w$ is fixed under each automorphism and ${\rm Aut}(F_n - v_{i_1}- \cdots - v_{i_s}) \cong {\rm Aut}(F_n - v_{i_1}- \cdots - v_{i_s}-w)$.

By using above point and regarding to the value of $t$, if we remove less than $t$ vertices of $F_n$, say $v_{i_1},\ldots , v_{i_{t-1}}$, then $D(F_n - v_{i_1}- \cdots - v_{i_{t-1}}) = D(F_n)$, and so $st_D(F_n) > t-1$. Therefore $D(F_n)=t$. \qed

With respect to the value of $st_D(F_n)$ for different values of $n$ we can prove the following result:
\begin{corollary}
For every natural number $k$, there exists a graph $G$ such that $st_D(G)=k$.
\end{corollary}
\proof If $k=1$, then it is suffice to consider the complete graphs. For $ k\geqslant 2$ the result follows immediately from Proposition \ref{pro fan}\qed

\section{Upper bounds for the distinguishing stability}
In this section, we study the relationship between the distinguishing stabilities of graphs $G$ and $G-v$, where $v\in V(G)$. Also we obtain upper bounds for $st_D(G)$. 
\begin{proposition}\label{minus}
Let $G$ be a graph and $v$ be a vertex of $G$, then $$st_D(G)\leqslant st_D(G-v)+1.$$
\end{proposition}
\proof   If $D(G)=D(G-v)$, then $st_D(G)\leqslant st_D(G-v)+1$. If $D(G)\neq D(G-v)$, then $st_D(G)= 1$, and so we have the result.\qed
 
 Using Proposition \ref{minus} and mathematical induction, we have $st_D(G)\leqslant st_D(G-v_1-\cdots - v_s)+s$ where $1\leqslant s \leqslant n-2$ and $n=\vert V(G)\vert$. Only using this formula we can get different upper bounds for $st_D(G)$ of graph $G$. We state some these upper bounds in the following theorem which for the prove of each case, it is sufficient to remove the vertices of $G$ until the induced subgraph which  stated in the hypothesis, appears. Next using $st_D(G)\leqslant st_D(G-v_1-\cdots - v_s)+s$ and the value of $st_D(G-v_1-\cdots - v_s)$, we can have the result. 
 
 \begin{theorem}
 Let $G$ be a simple graph of order $n\geqslant 2$.
 \begin{itemize}
 \item[(i)]  $st_D(G)\leqslant n-1$.
  \item[(ii)] If $\Delta$ is maximum degree of graph $G$, and $G$ has  the star graph $K_{1,\Delta}$ as the induced subgraph with $\Delta \geqslant 3$, then  $st_D(G)\leqslant n-\Delta$.
   \item[(iii)] If $d$ is the diameter of graph $G$, and $G$ has  the path $P_{d+1}$ as the induced subgraph, then  $$st_D(G)\leqslant \left\{ \begin{array}{ll}
   n-d+2 & d=3,\\
   n-d & d=1,4,\\
   n-d+1 & {\rm otherwise}.
   \end{array}\right.$$
  \item[(iv)]   $st_D(G)\leqslant n-{\rm \omega}(G)+1$, where $\omega(G)$ is the clique number of $G$.
 \end{itemize}
 \end{theorem}

\begin{theorem}
Let $G$ be a graph of order $n$, then $st_D(G)\leqslant n-D(G)+1$.
\end{theorem}
\proof 
Let $st_D(G)=k$. So for every $i$ ($1\leqslant i \leqslant k-1$) vertices of $G$, say $v_1,\ldots, v_{k-1}$,  we have $D(G)=D(G-v_1)=\cdots = D(G-v_1-\cdots - v_{k-1})$. Since the distinguishing number of a graph is at most equal to its order, so $D(G)=D(G-v_1-\cdots -v_{k-1})\leqslant n-k+1$. \qed

  A graph and its complement, always have the same automorphism group while their graph structure usually differs. Hence $D(G)=D(\overline{G})$ for every simple graph $G$. In the following theorem we use  this equality to show that the distinguishing stability of a graph and its complement are the same.
\begin{theorem}\label{complement}
Let $G$ be a simple graph, then $st_D(G)=st_D(\overline{G})$.
\end{theorem}
\proof 
Since $D(G)=D(\overline{G})$ and  $\overline{G}-v_1-\cdots - v_k=\overline{G -v_1 - \cdots - v_k}$, we have the result. \qed

\begin{theorem}\label{Dbound}
If there exists a vertex $v$ of $G$ such that $st_D(G)\leqslant st_D(G-v)$, then $st_D(G)\leqslant D(G)$.
\end{theorem}
\proof  If $D(G)\neq D(G-v)$, then $st_D(G)=1$, and so $st_D(G)\leqslant D(G)$. If $D(G)= D(G-v)$, then we use mathematical induction on the order of $G$. It can be seen that the result is true for small value of $n$. Let $st_D(G)\leqslant D(G)$ for all graphs of order $n < k$. Suppose that $n=k$, in this case by induction hypothesis we have $st_D(G-v)\leqslant D(G-v)$. Thus  we can write
\begin{equation*}
st_D(G) \leqslant st_D(G-v)\leqslant D(G-v) = D(G).
\end{equation*}

Therefore the result follows.\qed

\medskip 
We think that for any simple graph $G$ the inequality $st_D(G)\leqslant D(G)+1$ is true.  However,
until now all attempts to find a proof failed. So we propose the following conjecture here. 

\begin{conjecture}
Let $G$ be a simple connected graph, then $st_D(G)\leqslant D(G)+1$. 
\end{conjecture}

   Let $v$ be a vertex in $G$. The contraction of $v$ in $G$ denoted by $G/ v$ is the graph obtained by deleting $v$ and putting a clique on the (open) neighbourhood of $v$. Note that this operation does not create parallel edges; if two neighbours of $v$ are already adjacent, then they remain simply adjacent (see \cite{Walsh}). In the end of this section, we study  the distinguishing stabilities of graph $G/v$.  

\begin{theorem}
Let $G$ be a graph and $v$ be a vertex of it, and let $e_1,\ldots , e_k$ be the added edges to the neighbours of $v$ in the construction of $G/ v$. Suppose that $v_1,\ldots ,v_t$ are all different neighbours of $v$ which are incident to at least one of the edges  $e_1,\ldots , e_k$.  We have
\begin{equation*}
st_D(G/ v)\leqslant st_D(G-v-v_1-\cdots -v_t)+t.
\end{equation*}
\end{theorem}
\proof  First of all note that if $e$ is an edge incident to the vertex $w$ of $G$, then by Proposition \ref{minus}, $st_D(G-e)\leqslant st_D(G-w)+1$. Using Theorem \ref{complement} we can write
\begin{equation*}
st_D(G/ v)= st_D(G-v+e_1+\cdots +e_k)=st_D(\overline{G}-v-e_1-\cdots -e_k).
\end{equation*}

By the first sentence of proof, it can be concluded that 
\begin{equation*}
st_D(\overline{G}-v-e_1-\cdots -e_k)\leqslant st_D(\overline{G}-v-v_1-\cdots -v_t)+t.
\end{equation*}

Now using Theorem \ref{complement},  we obtain the result.\qed

\section{ Distinguishing bondage number of a graph }
In this section we study the edge stability (bondage) number on the distinguishing number of a graph.  The bondage number on the distinguishing number $b_D(G)$, (or the edge distinguishing  stability number) is the minimum number of edges whose removal changes the distinguishing number. First we determine $b_D(G)$ for several families of graphs including complete (bipartite) graphs, cycles, and paths.

\begin{proposition}
	\begin{enumerate}
		\item [(i)] For $n\geqslant 3$ and $m\geqslant 2$, we have $b_D(K_n)=b_D(K_{n,m})=1$. 
		
		\item[(ii)] 
	$b_D(P_n)=\left\{
		\begin{array}{ll}
		0& n=2\\
		2& n=3\\
		1& n=4\\
		2& n\geqslant 5.
		\end{array}\right. $
		\item[(iii)] 
		$b_D(C_n)=\left\{
		\begin{array}{ll}
			1 &n=3,4,5\\
			3& n\geqslant 6.
		\end{array}\right.
	$					\end{enumerate}
\end{proposition}
 
 \begin{theorem}
 For every natural number  $k$ and  $i\in \{1,\ldots , k+1\}$, there exists a graph $G_i$ such that $D(G_i)=2k$ and $b_D(G_i) =i$. Also, there exists a graph $H_i$ such that $D(H_i)=2k+1$ and $b_D(H_i) =i$.
 \end{theorem}
 \proof   With respect to Theorem \ref{distfan}, we know that for every $i\in \{1,\ldots , k+1\}$ there exists $n_i$ such that $D(F_{n_i})=2k$. Also it can be concluded that there exist $2k-1$ consecutive friendship graphs with the distinguishing number $2k$. Without loss of generality, we assume that $D(F_{n_i})= D(F_{n_i+1})=\cdots = D(F_{n_i+2k-2})=2k$. We claim that 
 \begin{itemize}
 \item[(i)] $b_D(F_{n_i+j})=j+1$ for $0\leqslant j \leqslant k-1$, and  
 \item[(ii)] $b_D(F_{n_i+k+j})=k+1$ for $0\leqslant j \leqslant k-2$.
 \end{itemize}  
 
Let $x_0$ be the central vertex and $x_{2t-1}$ and $x_{2t}$ ($t \geqslant 1$) be the two adjacent vertex on the base of $t$-th triangle of the friendship graph. 

 (i) If we remove the edges $x_0x_1,x_0x_3, \ldots , x_0x_{2j+1}$ from $F_{n_i +j}$, then $F_{n_i +j}$ converted to the friendship graph $F_{n_i-1}$ such that the end vertex of $j+1$ paths of order three identified to the central vertex of $F_{n_i-1}$. It can be computed that $D(F_{n_i +j}-x_0x_1 - x_0x_3 -  \cdots - x_0x_{2j+1})=2k-1$, and hence $b_D(F_{n_i +j})\leqslant j+1$. If we remove less than $j+1$ edges from $F_{n_i +j}$, say $e_1,\ldots , e_l$, $l\leqslant j$, then $F_{n_i +j} - e_1 - \cdots - e_l$ has $F_t$ where $t\geqslant n_i$ as its fixed induced subgraph, and hence $D(F_{n_i +j} - e_1 - \cdots - e_l)=2k$, so $b_D(F_{n_i +j})=j+1$.
 
 (ii) If we remove the edges $x_1x_2,x_3x_4, \ldots , x_{2k-1}x_{2k}$ from $F_{n_i+k +j}$, then we obtain a graph such that it has made by identifying the central vertices of the star graph $K_{1,2k}$ and  $F_{n_i+j-1}$, and so  $D(F_{n_i +k+j}-x_1x_2 - x_3x_4- \cdots - x_{2k-1}x_{2k})=2k-1$.  Hence $b_D(F_{n_i +k+j})\leqslant k+1$. If we remove less than $k+1$ edges from $F_{n_i+k +j}$, say $e_1,\ldots , e_l$, $l\leqslant k$, then $F_{n_i +k+j} - e_1 - \cdots - e_l$ has $F_{n_i+j+t}$ where $t\geqslant 0$ as its fixed induced subgraph, and hence $D(F_{n_i +j} - e_1 - \cdots - e_l)=2k$. Therefore  $b_D(F_{n_i +k +j})=k+1$.
 
 To prove the second part of theorem, it is sufficient to note that there exists $m_i$ such that $D(F_{m_i})= D(F_{m_i+1})=\cdots = D(F_{m_i+2k-1})=2k+1$. By a similar argument we can show that $b_D(F_{m_i+j})=j+1$, and $b_D(F_{m_i+k +j})=k+1$ for $0\leqslant j \leqslant k-1$. So we have the result.  \qed
 
 Similar to  Proposition \ref{minus}, we can prove $b_D(G)\leqslant b_D(G-e)+1$. Now we can state some upper bonds for $b_D(G)$.
 \begin{theorem}
 Let $G$ be a simple graph of order $n\geqslant 2$ and size $m$ with diameter $d$ and maximum degree $\Delta \geqslant 3$. Then we have 
 \begin{itemize}
 \item[(i)] $b_D(G) \leqslant m$.
  \item[(ii)] $b_D(G)\leqslant m-\Delta +1$.
   \item[(iii)] For $n\geqslant d+3$,  $b_D(G)\leqslant m-d+1$, and for $n < d+3$,  $b_D(G)\leqslant m-d+2$. 
  \item[(iv)] If  $\omega(G)$ is the clique number of $G$, then
  \begin{equation*}
b_D(G)\leqslant   \left\{
  \begin{array}{ll}
m-{\omega (G) \choose 2}+1  & {\rm if}~n\leqslant 2 \omega (G),\\
b_D(G)\leqslant m-{\omega (G) \choose 2}+ \omega (G)-1  & {\rm if}~n > 2 \omega (G).
  \end{array}\right.
  \end{equation*}
 \end{itemize}
 \end{theorem}
 \proof   Using $b_D(G)\leqslant b_D(G-e)+1$ and mathematical induction, we have $b_D(G)\leqslant b_D(G-e_1-\cdots - e_s)+s$ where $1\leqslant s \leqslant m-1$ and $m=\vert E(G)\vert$. To prove each case, it is sufficient to remove the edges of $G$ until the connected components of size greater than or equal  one has been stated in proof of each parts, appears. Next using $b_D(G)\leqslant b_D(G-e_1-\cdots - e_s)+s$ and the value of $b_D(G-e_1-\cdots - e_s)$, we can have the result.
 
 \begin{enumerate} 
 	\item[(i)] The   connected component of size greater than or equal one is the path $P_2$, so $b_D(G)\leqslant b_D(P_2\cup (n-2)K_1)+m-1$. Since $b_D(P_2\cup (n-2)K_1) =1$ so $b_D(G) \leqslant m$.

\item[(ii)]
 The   connected component of size greater than or equal  one is $K_{1,\Delta}$, so $b_D(G)\leqslant b_D(K_{1,\Delta}\cup (n-\Delta -1)K_1)+m-\Delta$. Since $b_D(K_{1,\Delta}\cup (n-\Delta -1)K_1)=1$, the result follows.

\item[(iii)]
 The   connected component of size greater than or equal one is $P_{d+1}$, so $b_D(G)\leqslant b_D(P_{d+1}\cup (n-d -1)K_1)+m-d$. If $n\geqslant d+3$, then $b_D(P_{d+1}\cup (n-d -1)K_1)=1$, and if  $n < d+3$, then $b_D(P_{d+1}\cup (n-d -1)K_1)\leqslant 2$. Hence the result follows.
 
 \item[(iv)]
  The   connected component of size greater than or equal one is $K_{\omega (G)}$, so $b_D(G)\leqslant b_D(K_{\omega (G)}\cup (n-\omega (G) )K_1)+m-{\omega (G) \choose 2}$. If $n\leqslant 2 \omega (G)$, then $b_D(K_{\omega (G)}\cup (n-\omega (G) )K_1)=1$, and if  $n > 2 \omega (G)$, then $b_D(K_{\omega (G)}\cup (n-\omega (G) )K_1)= \omega (G)-1$. Hence the result follows.  \qed
\end{enumerate}

In 1956, Nordhaus and Gaddum obtained the lower and upper
bounds for the sum of the chromatic numbers of a graph and its complement
(actually, the upper bound was first proved by Zykov \cite{Zykov} in 1949). Since then, Nordhaus-Gaddum type bounds were obtained for many graph
invariants. An exhaustive survey is given in \cite{Aouchiche}. Here, we  give Nordhaus-Gaddum type inequalities for the bondage number on the distinguishing number of a graph.

\begin{lemma}\label{Lemdifference}
If $G$ is a simple connected graph with an edge $e$ such that $b_D(G)\leqslant  b_D(G-e)$ and $b_D(\overline{G}+e)\leqslant  b_D(\overline{G})$, then $\vert b_D(G)-b_D(\overline{G})  \vert \leqslant 1$.
\end{lemma}
\proof The proof is by induction on $m$, the number of edges of $G$. It can easily be seen that it is  true for small value of $m$. So suppose that the result is  true 
for graphs with size less  than $m$. By $b_D(G)\leqslant b_D(G-e)+1$, we can obtain $\vert b_D(G)-b_D(\overline{G})  \vert \leqslant \vert b_D(G-e)-b_D(\overline{G}+e)  \vert$, and so using $b_D(\overline{G-e})=b_D(\overline{G}+e)$ and induction hypothesis we get the result.\qed
\begin{theorem}\label{Nordbandag}
If $G$ is a simple graph with an edge $e$ such that $b_D(G)\leqslant  b_D(G-e)$ and $b_D(\overline{G}+e)\leqslant  b_D(\overline{G})$, then $$1 \leqslant b_D(G)+b_D(\overline{G}) \leqslant 2 ~{\rm min}\{ b_D(G), b_D(\overline{G})\}+1.$$
\end{theorem}
\proof  This follows directly from Lemma \ref{Lemdifference}.\qed

Note that the upper bound of Theorem \ref{Nordbandag} is sharp and is achieved for example by the family of cycle graphs $C_n$ with $n\geqslant 7$. In fact, if $v_1,\ldots, v_n$ are consecutive vertices of $C_n$, then $\overline{C_n}-\{v_1,v_{n-1}\}-\{v_2,v_{n-3}\}$ is an asymmetric graph, and so $b_D(\overline{C_n})=2$. About the sharpness of lower bound, it is sufficient to consider the complete graphs.

\medskip
We conclude this section by mentioning a relation between the edge and vertex stability distinguishing number.

\begin{theorem}
Let $G$ be a graph. If there is an edge $e\in E(G)$  which is incident to one of the vertices of $G$ such that their removal change the distinguishing number and   $b_D(G)\leqslant b_D(G-e)$,  then $b_D(G) \leqslant \left\lceil \dfrac{st_D(G)}{2} \right\rceil +1$.
\end{theorem}
\proof  If there exists an edge $e$ such that $D(G) \neq D(G-e)$, then $b_D(G)=1$, and so we have the result. If $D(G)=D(G-e)$ for all edges of $G$, then by induction on the size of $G$, we can write
\begin{equation*}
b_D(G)\leqslant b_D(G-e)\leqslant \left\lceil \frac{st_D(G-e)}{2} \right\rceil +1 \leqslant \left\lceil \frac{st_D(G)}{2} \right\rceil +1.
\end{equation*}
where $e$ is incident to one of the vertices of $G$ such that their removal change the distinguishing number. Hence it can be concluded that $st_D(G-e)\leqslant st_D(G)$ for such edges. 
\qed
 
\end{document}